\documentclass[11pt]{article}

\usepackage{amssymb,amsmath}

\newtheorem{theorem}{Theorem}[section]

\newtheorem{Lemma}[theorem]{Lemma}
\newtheorem{defi}[theorem]{Definition}
\newtheorem{Corollary}[theorem]{Corollary}

\numberwithin{equation}{section}

\newcommand{\non}{\nonumber}
\newcommand{\cR}{\mathbb R}

\newcommand{\f}{\frac}

\newcommand{\eq}[1]{\mbox{\rm {(\ref{#1})}}}

\title{\Large\bf Regularity of solutions to a model for solid-solid phase transitions
driven by configurational forces}

\author{\small \sc Peicheng Zhu$^{1,2}$\thanks{E-mail: zhu@bcamath.org}\\
\small $^1$ Basque Center for Applied Mathematics (BCAM)\\
\small Building 500, Bizkaia Technology Park \\
\small E-48160 Derio, \ \ Spain\\
\small $^2$ IKERBASQUE, Basque Foundation for Science\\
\small E-48011 Bilbao,\ \ Spain
}

\date{ }

\begin{document}

\maketitle

\begin{abstract}

In a previous work, we prove the existence of weak solutions to an initial-boundary value
problem, with $H^1(\Omega)$ initial data, for a system of partial
differential equations, which consists of the
equations of linear elasticity and a nonlinear, degenerate parabolic
equation of second order. Assuming  in this article the initial data is in $H^2(\Omega)$, we investigate
 the regularity of weak solutions that is difficult due to the gradient term which plays a role
 of a weight. The problem models the behavior in time of materials
with martensitic phase transitions. This model with diffusive phase
interfaces was derived from a model with sharp interfaces, whose evolution is
driven by configurational forces, and can be thought to be a regularization
of that model. Our proof, in which the difficulties are caused by the weight in the
principle term, is only valid in one space dimension.



\end{abstract}
\section{Introduction}
Many inhomogeneous systems can be characterized by domains of different phases
separated by a distinct interface \cite{Emmerich03}. Driven out of
equilibrium, their dynamics result in the evolution of those interfaces which
might develop into structures (compositional and structural inhomogeneities) with characteristic length scales at the
nano-, micro- or meso-scale.  To a large extent, the material properties of such
systems are determined by those structures of small-scale. It is thus important to
understand precisely the mechanisms that drive the evolution of those structures.
 In this article we are interested in a model for the evolution, driven
 by configurational forces, of microstructures in  elastically deformable solids.
Materials microstructures may consist of spatially distributed phases of different
compositions and/or crystal structures, grains of different orientations, domains of
different structural variants, domains of different electrical or magnetic polarizations,
and structural defects, see e.g. \cite{Chen02}. These structural features usually have an intermediate
mesoscopic length scale in the range of nanometers to microns. The size, shape,
and spatial arrangement of the local structural features in a microstructure play a
critical role in determining the physical properties of a material. Because of the complex
and nonlinear nature of microstructure evolution, numerical approaches are often
employed.

There are two main types of modeling for the evolution of  microstructures.
 In the conventional approach, the regions separating the domains are treated as mathematically
sharp interfaces. The local interfacial velocity is then determined as part
of the boundary conditions, or is calculated from the driving force for interface motion
and the interfacial mobility. This approach requires the explicit tracking of the interface
positions. Such an interface-tracking approach can be successful in one-dimensional
systems, however it will be impractical for complicated three-dimensional microstructures.
Therefore, during the past decades, another approach
 has been invented, namely, the phase-field approach in which the interface is  not of zero
  thickness, instead an interfacial region with thickness of certain order of a small regularization
   parameter. Though it is still a young discipline in condensed matter physics,  this approach has emerged to be
 one of the most powerful methods for modeling the evolution of microstructures.
 It can be traced back  the theory of  diffuse-interface description, which is developed, independently, more
than a century ago by van der Waals \cite{VanderWaals} and some half century ago
by Cahn and Hilliard \cite{Cahn58}.

The two well-known models for temporal evolution of microstructures are
 the  Cahn-Hilliard/Allen-Cahn  equations corresponding, respectively, to the case that the order parameter is conserved and not
 conserved. These phase field models describe microstructure phenomena at the
mesoscale, and one suitable limit of it may be the corresponding sharp-  or thin-interface descriptions.
 In this article we study a model for the behavior in time of materials with
diffusionless phase transitions. The model has diffusive interfaces and
consists of the partial differential equations of linear elasticity coupled to
a quasilinear, non-uniformly parabolic equation of second order that differs from the Allen-Cahn equation
 (the Cahn-Hilliard equation in the case that the order parameter is conserved) by a
gradient term. It is derived in \cite{Alber00,Alber04} from a sharp interface model for diffusionless phase
transitions and can be considered to be a regularization of that model.
 To verify the validity of the new model,  mathematical analysis has been carried
 out for the existence of weak solutions to initial boundary value problems in
 one space dimension, \cite{Alber06,Alber07,Alber10a,Zhu}, the existence of spherically
 symmetric solutions \cite{Alber10b}, the motion of interfaces \cite{Alber09},
 and the existence of traveling waves \cite{Kawashima10}. In the present article, the existence and
 regularity
 of weak solutions to an initial-boundary value problem will be studied.
 We first formulate this initial-boundary value problem in the
three-dimensional case, reduce it
to the one-dimensional case and prove the existence of weak solutions to this one dimensional
problem, then study the regularity of weak solutions by assuming
that the initial data is in $H^2(\Omega)$.

\medskip
Let $\Omega\subset\cR^3 $ be an open set. It represents
the material points of a solid body. The different phases are
characterized by the order parameter $S(t,x)\in \cR$. A value of
$S(t,x)$ near to $ S_- $ (which is a real number) indicates that the material is in the matrix
phase at the point $x\in\Omega$ at time $t$, a value near to $ S_+ $
indicates that the material is in the second phase. The other unknowns
are the displacement $u(t,x)\in \cR^3 $ of the material point $x$ at
time $t$ and the Cauchy stress tensor $T(t,x)\in {\cal S}^3 $, where
${\cal S}^3 $ denotes the set of symmetric $3\times 3$-matrices. The
unknowns must satisfy the quasi-static equations
\begin{eqnarray}
 -{\rm div}_x\,T(t,x)&=&b(t,x),
 \label{eq1}  \\
 T(t,x) &=& D \big(\varepsilon(\nabla_x\, u(t,x)) -
            \bar\varepsilon S(t,x) \big),  \label{eq1a} \\
 S_t(t,x) &=& -c \Big( \psi_S(\varepsilon(\nabla_x
            u(t,x)),S(t,x)) - \nu \Delta_x S(t,x) \Big)|\nabla_x S(t,x)| \quad
 \label{eq2}
\end{eqnarray}
for $(t,x)  \in (0,\infty)\times\Omega$. The boundary and initial
conditions are
\begin{eqnarray}
 u(t,x)&=&\gamma(t,x),\ S(t,x)=0, \  (t,x)\in
 [0,\infty)\times\partial\Omega,
 \label{eq3a} \\
  S(0,x) &=& S_0(x),\  x\in \Omega.
 \label{eq3}
\end{eqnarray}
Here $\nabla_xu$ denotes the $3\times 3$-matrix of first order
derivatives of $u$, the deformation gradient, $(\nabla_xu)^T$ denotes the
transposed matrix and
$$
 \varepsilon(\nabla_xu) = \frac12\left(\nabla_x u+(\nabla_x u)^T\right)
$$
is the strain tensor. $\bar\varepsilon\in {\cal S}^3 $ is a given matrix,
the misfit strain, and $D:{\cal S}^3\to {\cal S}^3$ is
a linear, symmetric, positive definite mapping, the elasticity
tensor. In the free energy
\begin{equation}\label{freengy}
 \psi(\varepsilon,S)=\f12 \big(D(\varepsilon-\bar\varepsilon S)\big)
 \cdot (\varepsilon-\bar\varepsilon S)+\hat\psi(S)
\end{equation}
we choose for $\hat\psi\in C^2(\cR,[0,\infty))$ a double well potential with
minima at $S=S_-$ and $S=S_+$. The scalar product of two matrices is $A \cdot B =
\sum a_{ij} b_{ij}$. Also, $\psi_S$ is the partial derivative, $c>0$ is a
constant and $\nu$ is a small positive constant. Given are the volume force
$b:[0,\infty)\times\Omega\to \cR^3$ and the data $\gamma:[0,\infty)\times
\partial\Omega\to \cR^3$, $S_0: \Omega\to \cR$.

This completes the formulation of the initial-boundary value problem. The
equations (\ref{eq1}) and (\ref{eq1a}) differ from the system of linear
elasticity only by the term $\bar\varepsilon S$. The evolution equation
(\ref{eq2}) for the order parameter $S$ is non-uniformly parabolic because of
the term $\nu\Delta S |\nabla_xS|$. Since this initial-boundary value problem
is derived from a sharp interface model, to verify that it is indeed a
diffusive interface model regularizing the sharp interface model, it must be
shown that equations (\ref{eq1}) -- (\ref{eq3}) with positive $\nu$ have
solutions which exist global in time and is more regular if the initial data
 is more regular, and that these solutions tend to
solutions of the sharp interface model for $\nu\to 0$. This would also be a
method to prove existence of solutions to the original sharp interface model.

In this article we show that in one space dimension the initial-boundary value problem has
solutions, and study the regularity of these weak solutions with $H^2(\Omega)$ initial data.
 Whether solutions in three space dimensions exist and whether these
solutions converge to a solution of the sharp interface model for $\nu
\to 0$ is still an open problem to be investigated later. The model and
therefore the existence result is of interest not only in three dimensions but
also in one space dimension. Moreover we believe that this one-dimensional
existence result can also be helpful for an existence proof for
higher space dimensions.

Related to our investigations is the model for diffusion dominated phase
transformations obtained by coupling the elasticity equations (\ref{eq1}),
(\ref{eq1a}) with the Cahn-Hilliard equation. This model has recently
been studied in \cite{BCDGSS02, CMP00,Garcke03}.

\vskip0.2cm
\noindent {\bf Statement of the main result.} We
now assume that all functions only depend on the variables $x_1$ and $t$, and,
to simplify the notation, denote $x_1$ by $x$. The set $\Omega=(a,d)$ is a
bounded open interval with constants $a<d$. We write $Q_{T_e}:=(0,T_e) \times
\Omega$, where $T_e$ is a positive constant, and define
$$
 (v,\varphi)_Z = \int_Z v(y)\varphi(y)\, dy\,,
$$
for $Z = \Omega$ or $Z = Q_{T_e}$. If $v$ is a function defined on
$Q_{T_e}$ we denote the mapping $x \mapsto v(t,x)$ by $v(t)$. If
no confusion is possible we sometimes drop the argument $t$ and
write $v = v(t)$. We still allow that the material points can be
displaced in three directions, hence $u(t,x)\in\cR^3$, $T(t,x)\in
{\cal S}^3$ and $S\in\cR$. If we denote the first column of the matrix
$T(t,x)$ by $T_1(t,x)$ and set
$$
 \varepsilon(u_x)=\frac12\left((u_x,0,0)+(u_x,0,0)^T\right)\in {\cal S}^3,
$$
then with these definitions the equations (\ref{eq1}) -- (\ref{eq2}) in
the case of one space dimension can be written in the form
\begin{eqnarray}
 -T_{1x}&=&b,
 \label{Eeq1a1}\\
 T&=&D(\varepsilon(u_x)-\bar\varepsilon S),
 \label{Eeq1aa1}\\
 S_t&=&c\left(T\cdot\bar\varepsilon-\hat\psi^\prime(S)
 +\nu S_{xx}\right)|S_x|,
 \label{Eeq2a1}
\end{eqnarray}
which must be satisfied in $ Q_{T_e}$. Here we have inserted
 $\psi_S(\varepsilon, S) = -T\cdot \bar\varepsilon+ \hat\psi^\prime(S)$.
Since the equations (\ref{Eeq1a1}), (\ref{Eeq1aa1}) are linear, the
inhomogeneous Dirichlet boundary condition for $u$ can be reduced in the
standard way to the homogeneous condition. For simplicity we thus assume that
$\gamma=0$. The initial and boundary conditions therefore are
\begin{eqnarray}
 u(t,x)=0, && (t,x)\in (0,T_e)\times\partial\Omega,
 \label{Eeq3aa1}\\
 S(t,x)=0, && (t,x)\in (0,T_e)\times\partial\Omega,
 \label{Eeq3a1}\\
 S(0,x)=S_0(x), && x\in  \Omega.
 \label{Eeq4a1}
\end{eqnarray}

To define weak solutions of this initial-boundary value problem we note that
because of $\frac12(|y|y)^\prime=|y|$ equation (\ref{Eeq2a1}) is equivalent to
\begin{equation}
 S_t-c\nu\frac12 (|S_x|S_x)_x-c\left( T\cdot\bar\varepsilon-\hat\psi^\prime(S)
 \right)|S_x|=0.
\end{equation}

\begin{defi} \label{D1.1} Let $b\in L^\infty (0,T_e,L^2(\Omega))$,
$S_0\in L^\infty(\Omega)$. A function $(u, T, S)$ with
\begin{eqnarray}
 && u\in L^\infty(0,T_e;W^{1,\infty}_0(\Omega)) ,
 \label{property01} \\
 &&T\in L^\infty(Q_{T_e}),
 \label{property01a} \\
 && S\in L^\infty(Q_{T_e})\cap L^\infty(0,T_e; H^1_0(\Omega)),
 \label{property02}
\end{eqnarray}
 is a weak solution to  problem (\ref{Eeq1a1}) -- (\ref{Eeq4a1}), if  
 equations (\ref{Eeq1a1}) -- (\ref{Eeq1aa1}) with
(\ref{Eeq3aa1}) are satisfied in the weak sense, and if for all $\varphi\in
C^\infty_0((-\infty,T_e)\times \Omega)$
\begin{equation}
 (S,\varphi_t)_{Q_{T_e}}-c\nu\frac12 ( |S_x|S_x,\,\varphi_x )_{Q_{T_e}} +
  c\left( \big(T\cdot\overline{\varepsilon} - \hat{\psi}'(S) \big) |S_x|,\,
 \varphi\right)_{Q_{T_e}} + (S_0,\varphi(0))_\Omega = 0.
 \label{definition}
\end{equation}

\end{defi}

\medskip
We make the following assumption on the nonlinearity $\hat\psi(S)$.\\
{\bf Assumption for} $\hat\psi(S)$:
\begin{eqnarray}
 && \hat\psi(S) \in C^\infty(\cR)   \mbox{ is  a
 double-well  potential  which  has  two  local  minima at}\ \non\\
 && S_-  \mbox{  and}\ S_+   \mbox{  with}\ S_-<S_+ \mbox{  and one
 local maximum   at}\ S_* \mbox{  with}\ S_-<S_*<S_+ ,\quad \quad \    \label{B}\\
 &&  \mbox{  and  satisfies}\ \hat\psi'(S) >0 \mbox{  for}\ S_-< S < S_* \mbox{
 and}\ \hat\psi'(S) )<0  \mbox{  for}\
 S_*< S < S_+ .  \non
\end{eqnarray}
One typical example is: $\hat\psi(S) = (S(1-S))^2$ with $S_- = 0,\ S_+ = 1$.

\medskip
Now we are in a position to state the main result of this article.
%
\begin{theorem}\label{T1.1} Assume that $\hat\psi(S)$ satisfies \eq{B}.  Then for
 all $S_0\in H^1_0(\Omega)\cap H^2(\Omega)$ and $b \in
C(\overline{Q}_{T_e})$ with $b_t \in  C(\overline{Q}_{T_e})$, there
exists a weak solution $(u,T,S)$ to  problem (\ref{Eeq1a1}) --
(\ref{Eeq4a1}), which in addition to (\ref{property01}) --
(\ref{definition}) satisfies
\begin{equation}
 u_t\in C([0,T_e];H^1(\Omega)),\quad T_t\in C([0,T_e];L^2(\Omega))
 \label{proper1}
\end{equation}
and
\begin{equation}
 S_t\in L^\infty(0,{T_e};L^2(\Omega)),\quad
  |S_x|S_x \in L^\infty(0,{T_e}; H^1(\Omega)),\quad S_x\in  L^\infty(Q_{T_e}).
 \label{proper2}
\end{equation}

\end{theorem}

\bigskip
The remaining sections are devoted to the proof of this theorem. The main
difficulty in the proof stems from the fact that the coefficient $\nu|S_x|$ of
the highest order derivative $S_{xx}$ in equation (\ref{Eeq2a1}) is not
bounded away from zero and that it is not differentiable with respect to
$S_x$.

To prove Theorem~\ref{T1.1} we therefore consider in Section~2 a modified
initial-boundary value problem which consists of (\ref{Eeq1a1}),
(\ref{Eeq1aa1}), (\ref{Eeq3aa1}) -- (\ref{Eeq4a1}) and the equation
\begin{equation}
 S_t - c\nu |S_x|_\kappa S_{xx} - c
 \left(T\cdot\bar\varepsilon-\hat\psi^\prime(S) \right) (|S_x|_\kappa - \kappa)  =0 ,
 \  x\in \Omega,\  t>0
 \label{eq2a}
\end{equation}
with a constant $\kappa>0$. Here we use the notation
\begin{equation}
 |p|_\kappa:= {\sqrt{\kappa^2 + |p|^2}}.
 \label{pfunction}
\end{equation}
Since (\ref{eq2a}) is a uniformly (for $|S_x|\le M$,  $M$ is a positive constant)
 parabolic equation we can use a standard
theorem to conclude that the modified initial-boundary value problem has a
sufficiently smooth solution $(u^\kappa,T^\kappa,S^\kappa)$. For this solution
we derive in Section~3 a priori estimates independent of $\kappa\in (0,1]$.

The function $|p|$ is smoothed by $|p|_\kappa$ in \eq{pfunction} which is different from that in \cite{Alber06},
 thus we also prove the existence of weak solutions, though our main concern of this article
 is the regularity of solutions. To select a subsequence converging to a solution for $\kappa\to 0$ we need a
compactness result. However, our a priori estimates of $S_{xx}$ depend on a weight
$|S_x|_\kappa$, and  are not strong enough to
show that the sequence $S^\kappa_x$ is compact; instead, we can only show that
the sequence $\int_0^{S^\kappa_x}|y|^\frac12 dy$ is compact, from which we conclude
a subsequence of $S^\kappa_x$ that converges almost everywhere, thereby prove
the existence.   For the compactness proof in
 Section~4 we apply the compact Sobolev imbedding theorem, and don't need
  the Aubin-Lions Lemma or its  generalized form of this lemma
 given by Roub\'ic\v{e}k \cite{Roubicek}, Simon \cite{Simon}, which plays a crucial
  role in the article \cite{Alber06}.

In the proof of regularity estimates, we differentiate equation \eq{eq2a}
with respect to $t$. Thus a term like $(|S_x|_\kappa)_t \psi_S$ appears, which cannot
be absorbed by the a priori estimates with a weight $|S_x|_\kappa$. To overcome this
difficulty, we derive a type of estimate (see \eq{stsx}) with a reciprocal weight $|S_x|_\kappa^{-1}$.
 This is possible due to the special structure of the model studied here. However the Allen-Cahn model
 does not possesses such a structure, and our technique does not work for that model.

The method of proof is limited to one space dimension, since for the a priori
estimates it is crucial that the term $|S_x|S_{xx}$ in (\ref{Eeq2a1}) can be
written in the divergence form $\frac12(|S_x|S_{x})_x$. In the higher dimensional case
the corresponding term $|\nabla_xS|\Delta_xS$ cannot be rewritten in this
way, whence the multi-dimensional problem is still open.

\section{Existence of solutions to the modified problem}
In this section, we study the modified initial-boundary value problem and
show that it has a H\"older continuous classical solution. To formulate this
problem, let
$
\chi\in C_0^\infty(\cR,[0,\infty))
$
satisfy
$
\int_{-\infty}^{\infty} \chi(t)dt=1.
$
For $\kappa>0$, we set
$$
\chi_\kappa(t):=\frac1\kappa\chi\left(\frac{t}{\kappa} \right),
$$
and for $S\in L^\infty(Q_{T_e},\cR)$ we define
\begin{eqnarray}
 (\chi_\kappa*S)(t,x)=\int_{0}^{T_e} \chi_\kappa(t-s)S(s,x) ds.
 \label{convolution}
\end{eqnarray}
The modified initial-boundary value problem consists of the equations
\begin{eqnarray}
 -T_{1x}&=&b,
 \label{m2.1}\\
 T&=&D(\varepsilon(u_x)-\bar \varepsilon\chi_\kappa*S),
 \label{m2.2}\\
 S_t&=& c\nu|S_x|_\kappa S_{xx} + c(T\cdot \bar \varepsilon -
 \hat\psi^\prime(S))( |S_x|_\kappa - \kappa),
 \label{m2.3}
\end{eqnarray}
which must hold in $Q_{T_e}$, and of the boundary and initial conditions
\begin{eqnarray}
 u(t,x)&=&0,\quad (t,x)\in (0,T_e)\times \partial\Omega,
 \label{m2.4}\\
 S(t,x)&=&0, \quad (t,x)\in (0,T_e)\times \partial\Omega,
 \label{m2.5}\\
 S(0,x)&=&S_0(x), \quad x\in \Omega.
 \label{m2.6}
\end{eqnarray}

To formulate an existence theorem for this problem we need some
function spaces: For nonnegative integers $m,n$ and a real number
$\alpha\in (0,1)$ we denote by $C^{m+\alpha}(\overline{\Omega})$  the space of
$m-$times differentiable functions on $\overline{\Omega}$, whose $m-$th
derivative is H\"older continuous with exponent $\alpha$.  The space
$C^{\alpha,\alpha/2}(\overline{Q}_{T_e})$ consists of all functions on
$\overline{Q}_{T_e}$, which are H\"older continuous in the parabolic distance
$$
 {\rm d}((t,x),(s,y)):=\sqrt{|t-s|+|x-y|^2}.
$$
$C^{m,n}(\overline{Q}_{T_e})$ and $C^{m+\alpha,n+\alpha/2}(\overline{Q}_{T_e})$, respectively, are the spaces of
functions, whose $x$--derivatives up to order $m$ and $t$--derivatives up to
order $n$ belong to $C(\overline{Q}_{T_e})$ or to $C^{\alpha,\alpha/2}(\overline{Q}_{T_e})$, respectively.
%
%
\begin{theorem}\label{T2.1} Let $\nu,\kappa>0$, $T_e > 0$. Suppose
that the function $b \in C(\overline{Q}_{T_e})$ has the derivative
$b_t \in C(\overline{Q}_{T_e})$ and that the initial data $S_0\in
C^{2+\alpha}(\overline{\Omega})$ satisfy $S_0|_{\partial\Omega} =
S_{0,x}|_{\partial\Omega} = S_{0,xx}|_{\partial\Omega} = 0$. Then
there is a  solution
$$
 (u,T,S) \in C^{2,1}(\overline{Q}_{T_e}) \times
 C^{1,1}(\overline{Q}_{T_e}) \times
 C^{2+\alpha,1+\alpha/2}(\overline{Q}_{T_e})
$$
of the modified initial-boundary value problem (\ref{m2.1}) --
(\ref{m2.6}). This solution satisfies $S_{tx}\in L^2(Q_{T_e})$ and
\begin{eqnarray}
 \max_{\overline{Q}_{T_e}}|S|\le \max_{\overline{\Omega}}|S_0|.
 \label{m2.7}
\end{eqnarray}

\end{theorem}
\vskip0.2cm
{\it Proof.}   In \cite{Alber04} it is shown that the unique solution to the
 linear elliptic problem (\ref{m2.1}) -- (\ref{m2.2}), with
  (\ref{m2.4}) and given $S$,  is given by
\begin{eqnarray}
 u(t,x) &=& u^*\left( \int_a^x(\chi_\kappa*S)(t,y)dy-\frac{x-a}{d-a}
   \int_a^d(\chi_\kappa*S)(t,y)dy\right)+w(t,x),\qquad\mbox{}
   \label{m2.8a}\\
 T(t,x) &=& D(\varepsilon^*-\bar\varepsilon)(\chi_\kappa*S)(t,x)-
   \frac{D\varepsilon^*}{d-a} \int_a^d(\chi_\kappa*S)(t,y)dy + \sigma(t,x),
   \label{m2.8}
\end{eqnarray}
where $u^*\in\cR^3,\varepsilon^*\in {\cal S}^3$ are suitable constants only
depending on $\bar\varepsilon$ and $D$, and where for every $t\in [0,T_e]$
 the function $(w(t),\sigma(t)):\Omega\to \cR^3 \times {\cal S}^3$ is the
 solution to the boundary value problem
\begin{eqnarray*}
 -\sigma_{1x}(t) &=& b(t),  \\
 \sigma(t) &=& D\varepsilon(w_x(t)),  \\
 w(t)_{|\partial\Omega} &=& 0.
\end{eqnarray*}
Since by assumption $b$ and $b_t$ belong to $C(\bar Q_{T_e})$, it follows
 that $(w,\sigma)\in C^{2,1}(\bar Q_{T_e})
\times C^{1,1}(\bar Q_{T_e})$. We insert
(\ref{m2.8}) into (\ref{m2.3}) and obtain the equation
\begin{equation}
 S_t = a_1(S_x)S_{xx} + a_2\left(t,x,S,S_x,\chi_\kappa*S,
 \frac{1}{d-a} \int_a^d(\chi_\kappa*S)(t,y)dy \right)
 \label{m2.9}
\end{equation}
in $Q_{T_e}$\,, where
\begin{eqnarray*}
 a_1(p)=c\nu|p|_\kappa
\end{eqnarray*}
and
\begin{eqnarray*}
 a_2(t,x,S,p,r,s) = c\left(\bar\varepsilon\cdot D(\varepsilon^* -
  \bar\varepsilon)r - \bar\varepsilon \cdot D\varepsilon^* s
  + \bar\varepsilon\cdot \sigma(t,x)-\hat{\psi}'(S)\right)(|p|_\kappa - \kappa).
\end{eqnarray*}
Then (\ref{m2.9}), (\ref{m2.5}) and (\ref{m2.6})
form an initial-boundary value problem with nonlocal terms, which is
equivalent to the problem (\ref{m2.1}) -- (\ref{m2.6}).
We can apply \cite[Theorem~2.9, p.23]{Ladyzenskaya}, with a slight modification, to (\ref{m2.9}) and conclude
the existence of classical solution to the modified problem and   estimate (\ref{m2.7}) holds. For the details,
we refer to \cite{Alber06}.

%
\section{A priori estimates}

In this section we establish a-priori estimates for solutions of the modified
problem, which are uniform with respect to $\kappa\in (0,1]$. We remark that the
estimates in Lemma~\ref{L3.1} and some in Corollary~{C3.2}, though stated in the
one-dimensional case, can be generalized to higher space dimensions.

In what follows we assume that
\begin{equation}
 0<\kappa\le 1,
 \label{assumkappa}
\end{equation}
since we consider the limit $\kappa\to 0$. The $L^2(\Omega)$-norm is denoted
by $\|\cdot\|$, and the letter $C$ stands for varies positive constants
independent of $\kappa$, while may depend on $\nu$.

We start by constructing a family of approximate solutions to the
modified problem. To this end let $T_e$ be a fixed positive
number and choose for every $\kappa$ a function $S_0^\kappa\in
C_0^\infty(\Omega)$ such that
\begin{equation}
 \|S_0^\kappa-S_0\|_{H^1_0(\Omega)\cap H^2(\Omega)}\to 0,\quad \kappa\to 0,
 \label{approximate}
\end{equation}
where $S_0\in H^1_0(\Omega)\cap H^2(\Omega)$ are the initial data given in Theorem
\ref{T1.1}. We insert for $S_0$ in (\ref{m2.6}) the function $S_0^\kappa$ and
choose for $b$ in (\ref{m2.1}) the function given in
Theorem~1.1. These functions satisfy the assumptions of Theorem~\ref{T2.1},
hence there is a solution $(u^\kappa,T^\kappa,S^\kappa)$ of the
modified problem (\ref{m2.1}) -- (\ref{m2.6}), which exists in
$Q_{T_e}$. The inequality (\ref{m2.7}) and Sobolev's imbedding theorem
yield for this solution
\begin{equation}\label{approximateA}
 \sup_{0<\kappa\le 1}\|S^\kappa \|_{L^\infty(Q_{T_e})} \leq
 \sup_{0<\kappa\le 1} \|S_0^\kappa\|_{L^\infty(\Omega)} \le C.
\end{equation}
Remembering that $\sigma$ in (\ref{m2.8}) belongs to $C^{1,1}(\bar
Q_{T_e})$, we conclude from (\ref{approximateA}) that also
\begin{equation}
 \label{5max}
 \max_{\overline{Q}_{T_e}} |c(T^\kappa \cdot \overline{\varepsilon} -
 \hat{\psi}'(S^\kappa))| \leq C .
\end{equation}
%
\begin{Lemma} \label{L3.1} There hold  for any $t\in [0,T_e]$
\begin{eqnarray}
 \|S^\kappa_x(t)\|^2 + c\nu \int_0^t\int_\Omega
  |S^\kappa_x|_\kappa |S^\kappa_{xx}|^2dxd\tau &\le& C,
 \label{5}\\
 \int_0^t\int_\Omega \frac{|S^\kappa_t|^2}{|S^\kappa_x|_\kappa} dxd\tau &\le& C.
 \label{stsx}
\end{eqnarray}

\end{Lemma}

It is worth a remark on the estimate \eq{stsx}.\\
\noindent{\bf Remark 1.} {\it The reciprocal weight estimate \eq{stsx}: It is obtained by multiplying equation \eq{m2.3}
by the reciprocal of the weight function $|S^\kappa_x|_\kappa$ in the leading term. This is due to the special
structure of our model which makes it possible to get the regularity results in this article. In contract, the
model of the Allen-Cahn type with the mobility depending on the gradient of unknown
  does not possesses this feature, for instance, $S_t - |S_x|_\kappa S_{xx} + f(S) = 0$.
Suppose we obtain the $L^\infty(Q_{T_e})$-norm of $S$. If we multiply it by $S_t|S_x|^{-1}_\kappa $ and integrate it
with respect to $x$, we then see that the term
$\int_\Omega f(S) S_t|S_x|^{-1}_\kappa d x $ cannot be absorbed by the left hand side.
Consequently, it is impossible to obtain an estimate of
the form  \eq{stsx} for the Allen-Cahn type equation.

}

\medskip
\noindent{\it Proof.} Invoking $S^\kappa_{tx} \in L^2(Q_{T_e})$, by Theorem~\ref{T2.1},
which yields that for almost all $t$
$$
 \frac12\frac{d}{dt}\|S^\kappa_{x}(t)\|^2 = \int_\Omega S^\kappa_x(t)
 S^\kappa_{xt}(t) dx.
$$
Making use of this relation and estimate (\ref{5max}), multiplying (\ref{m2.3})
by $-S^\kappa_{xx}$ and integrating it with respect to $x$,  and
taking the boundary condition (\ref{m2.5}) into account, we obtain that
for almost all $t$
\begin{eqnarray}
 & &\frac12\frac{d}{dt}\|S^\kappa_{x}\|^2+
  c\nu \int_\Omega |S^\kappa_{x}|_\kappa |S^\kappa_{xx}|^2dx
 =\int_\Omega c(\hat{\psi}'(S^\kappa)-T^\kappa \cdot
 \overline{\varepsilon}) (|S^\kappa_{x}|_\kappa - \kappa) S^\kappa_{xx} dx \non
 \\[1ex]
 & \le & C\int_\Omega (|S^\kappa_{x}|_\kappa + \kappa) |S^\kappa_{xx}|dx
 = C\int_\Omega |S^\kappa_{x}|_\kappa^\frac12
  |S^\kappa_{x}|_\kappa^\frac12| S^\kappa_{xx}| dx + C\int_\Omega   \kappa  |S^\kappa_{xx}|dx \non\\[1ex]
 & \le& \frac{c\nu}4\int_\Omega|S^\kappa_{x}|_\kappa|S^\kappa_{xx}|^2 dx + \frac{c\nu}4\kappa\|S^\kappa_{xx}\|^2
 + \frac{ C }{\nu} \int_\Omega (|S^\kappa_x|_\kappa)^2 dx + C_\nu.
 \label{4}
\end{eqnarray}
Splitting the second term on the left hand side of  \eq{4} into two equal terms and  subtracting the term
$\frac{c\nu}4\int_\Omega|S^\kappa_{x}|_\kappa|S^\kappa_{xx}|^2 dx$ and $\frac{c\nu}4\kappa\|S^\kappa_{xx}\|^2$ on
both sides of this inequality, and using Gronwall's Lemma we derive
(\ref{5}) from the resulting estimate, noting also (\ref{approximate}) and $\kappa\le  |S^\kappa_x|_\kappa$.

\medskip
To derive \eq{stsx}, we multiply (\ref{m2.3}) by $ {S^\kappa_{t} }{|S^\kappa_{x}|_\kappa^{-1}} $ and integrate the resulting
equation with respect to $x$ to get
\begin{eqnarray}
 0 &=& \int_\Omega \frac{(S^\kappa_{t} )^2}{|S^\kappa_{x}|_\kappa} dx -
 c \int_\Omega ( \nu S^\kappa_{xx} - \psi_S) S^\kappa_{t} dx +
 c \int_\Omega \frac{\kappa \psi_S}{|S^\kappa_{x}|_\kappa }
 S^\kappa_{t}dx\non\\
 &=& \int_\Omega \frac{(S^\kappa_{t} )^2}{|S^\kappa_{x}|_\kappa}dx
 +I_1+I_2.
 \label{4a}
\end{eqnarray}
Invoking the formula $\psi_S = - T\cdot \bar\varepsilon +
\hat\psi'(S)$ and the boundary conditions, and using integration by
parts we have
\begin{eqnarray}
 I_1 &=& c\, \frac{d}{dt} \int_\Omega \left( \frac{\nu }{2} |S^\kappa_{x}|^2 + \hat\psi(S^\kappa)\right)
 dx - c\, \int_\Omega T^\kappa\cdot \bar\varepsilon S^\kappa_{t} dx\non\\
 &=& c\, \frac{d}{dt} \int_\Omega \left( \frac{\nu }{2} |S^\kappa_{x}|^2 + \hat\psi(S^\kappa)\right)
 dx + J.
 \label{4b}
\end{eqnarray}
To deal with the term $J$, we use \eq{5max} to get
\begin{eqnarray}
 |J| &=&  c\left|\int_\Omega T^\kappa \cdot \bar\varepsilon\, |S^\kappa_{x}|^\frac12_\kappa\, \frac{S^\kappa_{t}}{|S^\kappa_{x}|^\frac12_\kappa} dx\right|\non\\
 &\le&  C\left(\int_\Omega  |S^\kappa_{x}| _\kappa dx \right)^\frac12
 \left(\int_\Omega  \frac{(S^\kappa_{t})^2}{|S^\kappa_{x}| _\kappa} dx\right)^\frac12\non\\
 &\le&  C\left(\|S^\kappa_{x}\| + 1 \right)^\frac12
 \left(\int_\Omega  \frac{(S^\kappa_{t})^2}{|S^\kappa_{x}| _\kappa} dx\right)^\frac12 \non\\
 &\le& \frac12  \int_\Omega  \frac{(S^\kappa_{t})^2}{|S^\kappa_{x}| _\kappa} dx + C .
 \label{4c}
\end{eqnarray}
Here we used   the estimate \eq{5} and the Cauchy-Schwarz and Young inequalities.

For $I_2$, we make use of equation \eq{m2.3} and write
\begin{eqnarray}
 I_2 &=& c \int_\Omega \frac{\kappa \psi_S}{|S^\kappa_{x}|_\kappa }
 \Big(c\nu |S^\kappa_{x}|_\kappa S^\kappa_{xx} - c\psi_S (|S^\kappa_{x}|_\kappa - \kappa)\Big) dx\non\\
 &=& c^2 \int_\Omega \left( \nu \kappa \psi_S
   S^\kappa_{xx} - \kappa (\psi_S)^2 \frac{ |S^\kappa_{x}|_\kappa - \kappa}{|S^\kappa_{x}|_\kappa}\right)
   dx.
 \label{4d}
\end{eqnarray}
By definition, one has $|S^\kappa_{x}|_\kappa\ge \kappa$ which
implies $\frac{\kappa}{|S^\kappa_{x}|_\kappa}\le 1$. So
\begin{eqnarray}
 |I_2| & \le & C \int_\Omega  (\kappa
   |S^\kappa_{xx}| +  |S^\kappa_{x}|_\kappa + \kappa ) dx\non\\
   & \le & c\nu \kappa \|S^\kappa_{xx}\|^2 + C_\nu(\|S^\kappa_{x}\| + 1 ) .
 \label{4e}
\end{eqnarray}

With the help of  \eq{5}, \eq{4b} -- \eq{4e}, we integrate \eq{4a} with respect
to $t$, then obtain
\begin{eqnarray}
 \frac12 \int_0^t\int_\Omega \frac{(S^\kappa_{t} )^2}{|S^\kappa_{x}|_\kappa}dx
 + c \int_\Omega \left(  \frac{\nu }{2} |S^\kappa_{x}|^2 + \hat\psi(S^\kappa) \right) dx\le C,
 \label{4f}
\end{eqnarray}
which implies \eq{stsx}. Thus we complete the proof of this lemma.

\medskip
Furthermore, we obtain
%
%
\begin{Corollary}\label{C3.2} There holds for any $t\in
[0,T_e]$
\begin{eqnarray}
 \int_0^t\int_\Omega\left(|S^\kappa_x|_\kappa|S^\kappa_{xx} |  \right)^\frac43dxd\tau &\le& C,
 \label{5d}\\
 \int_0^t\int_\Omega\left(|S^\kappa_x S^\kappa_{xx} |  \right)^\frac43dxd\tau &\le& C,
 \label{5w}\\
 \int_0^t\left\|\int_0^{S^\kappa_x}|y|_\kappa dy\right\|_{W^{1,\frac43}(\Omega)}^\frac43 d\tau &\le& C,
 \label{5g}\\
 \int_0^t \left\|\int_0^{S^\kappa_x}|y|_\kappa dy\right\|_{L^\infty(\Omega)}^\frac43 d\tau &\le& C,
 \label{5e}\\
 \|\, |S^\kappa_x| S^\kappa_x \|_{L^\frac43(0,T_e;L^\infty(\Omega))} &\le& C,
 \label{5h}\\
 \int_0^t \left\|S^\kappa_x\right\|_{L^\infty(\Omega)}^\frac83 d\tau &\le& C.
 \label{5f}
\end{eqnarray}

\end{Corollary}

\noindent{\it Proof.} For some $2 > p \ge 1$ we choose $q,\ q'$ such that
$$
 q=\frac2p, \quad \frac1q + \frac{1}{q^\prime}=1.
$$
By H\"older's inequality,  we have
\begin{eqnarray}
& & \int_0^t\int_\Omega\left(|S^\kappa_{x}|_\kappa|S^\kappa_{xx}|
  \right)^pdx d\tau \non\\
&= & \int_0^t\int_\Omega \left(|S^\kappa_{x}|_\kappa\right)^\frac{p}{2}
  \left( \left(|S^\kappa_{x}|_\kappa\right)^\frac{p}{2}
  |S^\kappa_{xx} |^p\right) dx d\tau \non \\
&\le&  \left(\int_0^t\int_\Omega
  \left(|S^\kappa_{x}|_\kappa \right)^\frac{pq^\prime}{2}  dxd\tau
  \right)^\frac{1}{q^\prime}\left(\int_0^t\int_\Omega
  \left(|S^\kappa_{x}|_\kappa \right)^\frac{pq}{2}
  |S^\kappa_{xx}|^{pq}  dxd\tau \right)^\frac{1}{q} \non\\
&\le &  \left(\int_0^t \int_\Omega
  \left(|S^\kappa_{x}|_\kappa\right)^\frac{p}{2-p}
  dx d\tau\right)^\frac{2-p}{2}
  \left(\int_0^t \int_\Omega |S^\kappa_{x}|_\kappa|S^\kappa_{xx} |^2
  dx d\tau \right)^\frac{p}{2}.
\label{6}
\end{eqnarray}
Estimate (\ref{5}) implies that if $p$ satisfies $\frac{p}{2-p}\le 2$, i.e. $p\le
 \frac43$,  then the right hand side of (\ref{6}) is bounded. This
yields  estimate \eq{5d}.
 By definition of $|y|_\kappa$,
$$
  |y|_\kappa - |y| = \frac{\kappa^2}{|y|_\kappa + |y|}.
$$
Since   $|y|_\kappa+|y|\ge \kappa $, we have
\begin{eqnarray}
 \frac{\kappa^2}{|y|_\kappa+|y|}  \le \frac{\kappa^2}{\kappa}= \kappa.
 \label{2.6a2}
\end{eqnarray}
Hence
$$
0\le |y|_\kappa - |y| \le \kappa.
$$
Letting $y = S^\kappa_{x}$ yields
$$
 |S^\kappa_{x} S^\kappa_{xx}| = |S^\kappa_{x}| |S^\kappa_{xx}| \le (|S^\kappa_{x}|_\kappa - |S^\kappa_{x}|) |S^\kappa_{xx}|\le \kappa |S^\kappa_{xx}|,
$$
and \eq{5w} follows from \eq{5d} and estimate \eq{5}.

Next we are going to prove  \eq{5g}. Writing
\begin{eqnarray}
 |S^\kappa_x|_\kappa S^\kappa_{xx}= \left(\int_0^{S^\kappa_x}|y|_\kappa dy\right)_x,
 \label{sxxsx}
\end{eqnarray}
 and invoking that the primitive of $|y|_\kappa$ is equal to
$$
\frac12\left(y\sqrt{y^2+ \kappa^2} + \kappa^2\log\Big(y + \sqrt{y^2+ \kappa^2} \Big)  \right),
$$
which,  thanks to $\log x\le  x -1 $ for all $x>0$, is bounded by $ C(y^2 + 1)$,  we then show easily that
$$
\int_\Omega \int_0^{S^\kappa_x}|y|_\kappa dy dx\le C\int_\Omega (|S^\kappa_x|^2 +1)dx\le C.
$$
To apply the Poincar\'e  inequality of the form
$$
 \|f-\bar f\|_{L^p(\Omega)}\le C \|f_x\|_{L^p(\Omega)}
$$
where  $\bar f:=\frac1{|\Omega|}\int_\Omega f(x)dx$, we choose
$$
 p=\frac43,\quad f = \int_0^{S^\kappa_x}|y|_\kappa dy ,
$$
and obtain
\begin{eqnarray}
 &&\int_0^t \left\|\int_0^{S^\kappa_x}|y|_\kappa dy  \right\|^\frac43_{L^{\frac43} (\Omega) }d\tau  \non\\
 &\le& C\int_0^t \left\| \left(\int_0^{S^\kappa_x}|y|_\kappa dy \right)_x \right\|^\frac43_{L^{\frac43} (\Omega) }d\tau
 + C \int_0^t \left\|\, {\overline{ \int_0^{S^\kappa_x} |y|_\kappa   dy} }\, \right\|^\frac43_{L^{\frac43} (\Omega) } d\tau \non\\
 &\le& C\int_0^t \left\| \, |S^\kappa_x |_\kappa S^\kappa_{xx} \right\|^\frac43_{L^{\frac43} (\Omega) } d\tau  + C \int_0^t  1\, d\tau,
 \label{sxIntegral}
\end{eqnarray}
which implies, by \eq{5d}, that
\begin{eqnarray}
 \int_0^t \left\|\int_0^{S^\kappa_x}|y|_\kappa dy  \right\|^\frac43_{L^{\frac43} (\Omega) }d\tau
 &\le& C.
  \label{sxIntegral2}
\end{eqnarray}
Hence \eq{5g} follows, and we get $\int_0^{S^\kappa_x}|y|_\kappa dy\in L^\frac43(0,T_e;W^{1,\frac43}(\Omega))$.
Making use of the Sobolev embedding theorem, we get \eq{5e}.

It remains to prove estimate (\ref{5f}), since (\ref{5h}) is
equivalent to (\ref{5f}).
 We rewrite $\int_0^{S^\kappa_x}|y|_\kappa dy$ as
\begin{eqnarray}
\int_0^{S^\kappa_x}|y|_\kappa dy
 &= &\int_0^{S^\kappa_x}|y| dy+ \int_0^{S^\kappa_x}(|y|_\kappa-|y|) dy\non\\
 &= &\left.\frac12|y| y\right|_0^{S^\kappa_x}
+ \int_0^{S^\kappa_x}\frac{\kappa^2}{|y|_\kappa+|y|} dy\non\\
 &= &\frac12 |S^\kappa_x| S^\kappa_x
+ \int_0^{S^\kappa_x}\frac{\kappa^2}{|y|_\kappa+|y|} dy.
\label{2.6a0}
\end{eqnarray}
Thus
\begin{eqnarray}
 \frac12 (|S^\kappa_x| S^\kappa_x)_x
 &= &\left(\int_0^{S^\kappa_x}|y| dy\right)_x =
 \left(\int_0^{S^\kappa_x}|y|_\kappa dy\right)_x
 - \frac{\kappa^2 S^\kappa_{xx}}{|S^\kappa_x|_\kappa+|S^\kappa_x|}.
 \label{2.6a1}
\end{eqnarray}
By \eq{2.6a2} and the Young inequality we obtain from (\ref{5}) and the assumption that $k\le 1$ that
\begin{eqnarray}
 \left|\frac{\kappa^2
 S^\kappa_{xx}}{|S^\kappa_x|_\kappa+|S^\kappa_x|}\right| &\le&
 |\kappa  S^\kappa_{xx}|,\ {\rm thus}\non\\
 \|\kappa S^\kappa_{xx}\|_{L^\frac43(Q_{T_e})} &\le&
 \left(\int_{Q_{T_e}} \left(\kappa^2
 +\kappa|{S^\kappa_{xx}}|^2\right)d xd\tau  \right)^\frac34\le C.
\end{eqnarray}
Combination with \eq{5g}  and  (\ref{2.6a1}) yields
\begin{eqnarray}
 \| (|S^\kappa_x| S^\kappa_x)_x\|_{L^\frac43(Q_{T_e})}
 \le C \left\|\left(\int_0^{S^\kappa_x}|y|_\kappa dy\right)_x\right\|_{L^\frac43(Q_{T_e})}
 + C \|\kappa S^\kappa_{xx}\|_{L^\frac43(Q_{T_e})}  \le C.
 \label{sxsxx}
\end{eqnarray}

It is clear that $ {\overline {|S^\kappa_x| S^\kappa_x} } \le C \int_\Omega |S^\kappa_x| ^2 dx\le C$. Applying
 again the Poincar\'e  inequality to the function $f=|S^\kappa_x| S^\kappa_x$, we arrive at
$$
 \left\|\, |S^\kappa_x| S^\kappa_x\right\|_{L^\frac43(Q_{T_e})} \le C.
$$
Hence this, combined with \eq{sxsxx}, implies that
$$
 \|\, |S^\kappa_x| S^\kappa_x\|_{L^\frac43(0,{T_e};W^{1,\frac43 }(\Omega))} \le C,
$$
one concludes by using the Sobolev embedding theorem that
$$
\|\, |S^\kappa_x| S^\kappa_x\|_{L^\frac43(0,{T_e};L^{\infty}(\Omega))} \le C,
$$
which is
$$
 \| S^\kappa_x\|_{L^\frac83(0,{T_e};L^{\infty}(\Omega))} \le C.
$$
This completes the proof of this corollary.
%

\begin{Lemma} \label{L3.3} There hold  for any $t\in [0,T_e]$
\begin{eqnarray}
 \|S^\kappa_t(t)\|^2 + c\nu \int_0^t\int_\Omega
   |S^\kappa_x|_\kappa |S^\kappa_{xt}|^2dxd\tau &\le& C,
 \label{stsxt}\\
 \|\, |S^\kappa_x|_\kappa S^\kappa_{xx}(t)\| &\le& C.
 \label{sxsxx2}
\end{eqnarray}

\end{Lemma}

\noindent{\it Proof.} Suppose that estimate \eq{stsxt} is true, from equation
\eq{m2.3} and the estimate \eq{5}   we can easily get \eq{sxsxx2}. So it is enough to prove
\eq{stsxt}.

Differentiating \eq{m2.3} formally with respect to $t$ yields
\begin{eqnarray}
 S^\kappa_{tt}&=& c\nu(|S^\kappa_x|_\kappa S^\kappa_{xt})_x
 + c\left((T^\kappa\cdot \bar \varepsilon -
 \hat\psi^\prime(S^\kappa))( |S^\kappa_x|_\kappa - \kappa)\right)_t.
 \label{stt}
\end{eqnarray}
Multiplying \eq{stt} by $S^\kappa_{t}$ and integrating it,  by
integration by parts, we obtain
\begin{eqnarray}
 0 &=& \frac12 \frac{d}{dt} \|S^\kappa_{t}\|^2
 + c\nu\int_\Omega |S^\kappa_x|_\kappa |S^\kappa_{xt} |^2dx
 + c\int_\Omega\left((T^\kappa\cdot \bar \varepsilon -
 \hat\psi^\prime(S^\kappa))( |S^\kappa_x|_\kappa - \kappa)\right)_t S^\kappa_{t}dx\non\\
 &=& \frac12 \frac{d}{dt} \|S^\kappa_{t}\|^2
 + c\nu\int_\Omega |S^\kappa_x|_\kappa |S^\kappa_{xt} |^2dx
 + J_1.
 \label{stt1}
\end{eqnarray}
It is not difficult to  carry out a rigorous justification of \eq{stt1} with the help of
difference quotient, we omit the details. Computation gives
\begin{eqnarray}
 J_1 &=& c\int_\Omega\left( (T^\kappa\cdot \bar \varepsilon -
 \hat\psi^\prime(S^\kappa)) _t( |S^\kappa_x|_\kappa - \kappa)
 + (T^\kappa\cdot \bar \varepsilon -
 \hat\psi^\prime(S^\kappa) )( |S^\kappa_x|_\kappa  )_t \right)
 S^\kappa_{t}dx\non\\
 &=& J_{11} + J_{12}.
 \label{stt2}
\end{eqnarray}
By the formula of $T$, we have
\begin{eqnarray}
 |J_{11}|
 &\le& C\int_\Omega\left(|S^\kappa_{t}|^2+|S^\kappa_{t}|\right)( |S^\kappa_x|_\kappa +
 \kappa)dx\non\\
 &\le& C( \|S^\kappa_x\|_{L^\infty(\Omega)} +
 1)\int_\Omega(|S^\kappa_{t}|^2+1) dx\non\\
 &\le& C( \|S^\kappa_x\|_{L^\infty(\Omega)} +
 1)(\|S^\kappa_{t}\|^2+1) .
 \label{stt3a}
\end{eqnarray}
To handle $J_{12}$, we make use of estimate \eq{stsx} and $|y| \le
|y|_\kappa$.
\begin{eqnarray}
 |J_{12}| &\le& C\int_\Omega
 \frac{|S^\kappa_x
 S^\kappa_{xt}S^\kappa_{t}|}{|S^\kappa_x|_\kappa}dx
 = C\int_\Omega
 \frac{|S^\kappa_x S^\kappa_{xt}|}{|S^\kappa_x|^\frac12_\kappa}\,
 \frac{|S^\kappa_{t}|}{|S^\kappa_x|^\frac12_\kappa}dx\non\\
 &\le& C\int_\Omega
 \frac{|S^\kappa_x|_\kappa
 |S^\kappa_{xt}|}{|S^\kappa_x|^\frac12_\kappa}\,
 \frac{|S^\kappa_{t}|}{|S^\kappa_x|^\frac12_\kappa}dx\non\\
 &\le&  \frac{c\nu}2 \int_\Omega
  |S^\kappa_x|_\kappa| S^\kappa_{xt} |^2 dx +
 C_\nu \int_\Omega\frac{|S^\kappa_{t}|^2}{|S^\kappa_x|_\kappa}dx.
 \label{stt3b}
\end{eqnarray}

Thus it follows from \eq{stt1} -- \eq{stt3b} that
\begin{eqnarray}
 & & \frac12 \frac{d}{dt} \|S^\kappa_{t}\|^2
 + c\nu\int_\Omega |S^\kappa_x|_\kappa |S^\kappa_{xt} |^2dx\non\\
 &\le & \frac{c\nu}2 \int_\Omega
  |S^\kappa_x |_\kappa|S^\kappa_{xt}|^2 dx +
 C_\nu \int_\Omega\frac{|S^\kappa_{t}|^2}{|S^\kappa_x|_\kappa}dx + C( \|S^\kappa_x\|_{L^\infty(\Omega)} +
 1)(\|S^\kappa_{t}\|^2+1).
 \label{stt4}
\end{eqnarray}
From equation \eq{m2.3} and assumption $S_0\in H^2(\Omega)$ we compute the initial data
\begin{eqnarray}
 \|S^\kappa_{t}|_{t=0}\| &\le& C(\|\,|S_{0x}|_\kappa S_{0xx} \| +
 \| |S_{0x}|_\kappa + \kappa \| )\non\\
 &\le &C(\|\,|S_{0x}|_\kappa\|_{L^\infty(\Omega)} \|  S_{0xx} \| +
 \|  S_{0x}\| + 1 )\non\\
 &\le &C((\|S_{0x }\|_{H^1(\Omega)} + 1) \|  S_{0xx} \| +
 \|  S_{0x}\| + 1 )\le  C.
 \label{stinitial}
\end{eqnarray}
Thus $S^\kappa_{t}|_{t=0}\in L^2(\Omega)$. Next we use the Gronwall inequality of the form:
\begin{Lemma} \label{Gronwall}
For measurable functions $y,A,B$ defined on $[0,T_e]$, such that $y\ge 0$ and $ A,B\in L^1(0,T_e)$, if
$$
 y'(t)\le A(t)y(t) + B(t),
$$
then
$$
 y(t)\le y(0)\exp\left(\int_0^t A(\tau) d\tau\right) + \int_0^t B(s) \exp\left( \int_s^t A(\tau) d\tau \right) ds.
$$

\end{Lemma}
Defining
$$
 y(t)=  \|S^\kappa_{t}(t)\|^2,\ A(t) = C(\|S^\kappa_x\|_{L^\infty(\Omega)} + 1),\  B(t) = C( \|S^\kappa_x\|_{L^\infty(\Omega)} +
 1) + C_\nu \int_\Omega\frac{|S^\kappa_{t}|^2}{|S^\kappa_x|_\kappa}dx,
$$
where $A,\ B$  are integrable over $[0,T_e]$ by Lemma~\ref{L3.1} and Corollary~\ref{C3.2}, we derive from \eq{stt4} and \eq{stinitial} that
\begin{eqnarray}
   \|S^\kappa_{t}(t)\|^2
 + c\nu\int_0^t\int_\Omega |S^\kappa_x|_\kappa |S^\kappa_{xt} |^2dxd\tau
 \le C \|S^\kappa_{t}(0)\|^2 + C  \le C .
 \label{stt5}
\end{eqnarray}
Thus the proof of this lemma is complete.

\begin{Corollary} \label{C3.5}
The function $\int_0^{S^\kappa_x}|y|_\kappa dy$ belongs to $H^1(Q_{T_e})$, and the estimates  hold
\begin{eqnarray}\label{state0}
 \left\| \left(\int_0^{S^\kappa_x}|y|^\frac12_\kappa dy\right)_t \right\|_{L^{2}(Q_{T_e})} \leq C\,, \\[1ex]
 \left\| \int_0^{S^\kappa_x}|y|^\frac12_\kappa dy \right\|_{L^2(0,T_e;H^{1}(\Omega))} \le C\,.
 \label{state1}
\end{eqnarray}

\end{Corollary}

\noindent{\it Proof.}  Calculating yields
\begin{eqnarray}
 \left(\int_0^{S^\kappa_x}|y|^\frac12_\kappa dy\right)_t = |S^\kappa_{x}|^\frac12_\kappa S^\kappa_{xt},
 \label{state2}
\end{eqnarray}
recalling \eq{stsxt}, we obtain \eq{state0}. Similarly,
\begin{eqnarray}
 \left(\int_0^{S^\kappa_x}|y|^\frac12_\kappa dy\right)_x = |S^\kappa_{x}|^\frac12_\kappa S^\kappa_{xx},
 \label{state2a}
\end{eqnarray}
combining with \eq{5} gives
\begin{eqnarray}
 \left\|\left(\int_0^{S^\kappa_x}|y|^\frac12_\kappa dy\right)_x \right\|_{L^2(Q_{T_e})} \le C .
 \label{state2b}
\end{eqnarray}
Finally,  Noting $|\int_0^{S^\kappa_x}|y|^\frac12_\kappa dy| \le C\max\left\{M,|S^\kappa_x|^\frac32\right\}$
for some large constant $M>0$, we have
\begin{eqnarray}
 \left\| \int_0^{S^\kappa_x}|y|^\frac12_\kappa dy  \right\|_{L^2(Q_{T_e})}^2
 & \le& C +  C\int_\Omega |S^\kappa_x|^3dx\non\\
 & \le& C +  C\|S^\kappa_x\|_{L^\infty(\Omega)}\int_\Omega |S^\kappa_x|^2dx .
 \label{state2c}
\end{eqnarray}
Thus by \eq{5f} in Corollary~\ref{C3.2} there holds
\begin{eqnarray}
 \int_0^t\left\| \int_0^{S^\kappa_x}|y|^\frac12_\kappa dy  \right\|_{L^2(Q_{T_e})}^2d\tau
 & \le&  C +  C\int_0^t\|S^\kappa_x\|_{L^\infty(\Omega)}d\tau \le C.
 \label{state2d}
\end{eqnarray}
Then \eq{state1} follows from \eq{state2b} and \eq{state2d}. The proof of the lemma is complete.

%
%
\section{Existence/regularity of solutions to the phase field model}

We shall make use of  the a priori estimates established in the previous
section to study the convergence of
$(u^\kappa,T^\kappa,S^\kappa)$ as $\kappa\to 0$. In this section we will
show that there is a subsequence, which converges to a weak solution of the
initial-boundary value problem (\ref{Eeq1a1}) -- (\ref{Eeq4a1}), thus we prove
the existence of weak solutions; then we shall investigate the regularity of solutions.

\medskip
\noindent{\bf Existence.} It follows from Lemmas~\ref{L3.1} and  \ref{L3.3} that
\begin{equation}\label{m4.1}
 \|S^\kappa\|_{H^{1}(Q_{T_e})} \leq C\,,
\end{equation}
for a constant $C$ independent of $\kappa$. Hence, we can select a sequence
$\kappa_n \rightarrow 0$ and a function $S\in H^{1}(Q_{T_e})$,
such that the sequence $S^{\kappa_n}$, which we again denote by $S^\kappa$,
satisfies
\begin{equation}
 \label{m4.2}
 \| S^\kappa - S \|_{L^{2}(Q_{T_e})} \rightarrow 0,\qquad
 S^\kappa_x \rightharpoonup S_x\,,\qquad S^\kappa_t \rightharpoonup S_t\,,
\end{equation}
where the weak convergence is in $L^{2}(Q_{T_e})$\,.

Since the equation (\ref{Eeq2a1}) is nonlinear, the weak convergence
of $S^\kappa_x$ is not enough to prove that the limit function solves this
equation. In the following lemma we therefore show that $S^\kappa_x$ converges
pointwise almost everywhere:

\begin{Lemma} \label{L4.1}
There exists a subsequence of $S_x^\kappa$,
we still denote it by $S_x^\kappa$,
such that
\begin{eqnarray}
 \label{m4.3}
 \int_0^{S^\kappa_x}|y|^\frac12_\kappa dy \to
  \int_0^{S_x}|y|^\frac12  dy & & {  a.e. \ \   in\  \ } Q_{T_e},\\
 \label{m4.4}
 S_x^\kappa\to S_x,   & & a.e. \  \   in\  \  Q_{T_e},\\[0.2cm]
 \label{m4.5}
 |S_x^\kappa|_\kappa\rightharpoonup
 |S_x|,   & &
 {weakly \ in }\  L^2(Q_{T_e}), \\
 \label{m4.3a}
 \int_0^{S^\kappa_x}|y|^\frac12 dy \to \int_0^{S_x}|y|^\frac12  dy, & &  strongly
 \  \  in \  \ {L^2(Q_{T_e} )},
\end{eqnarray}
 as $\kappa\to 0$.
\end{Lemma}

The proof is based on  the following  result:
\begin{Lemma} \label{L4.2}
Let $(0,T_e)\times \Omega$  be an open
set in $\cR^+\times \cR^n$.  Suppose
functions $g_n, g$ are in $L^q((0,T_e)\times \Omega )$ for any given  $1<q
<\infty$, which satisfy
$$
 \|g_n\|_{L^q((0,T_e)\times \Omega )}\le C, \ \ g_n\to g \ almost\ everywhere \
 in\ (0,T_e)\times \Omega .
$$
Then $g_n$ converges to $g$ weakly in $L^q((0,T_e)\times \Omega )$.
\end{Lemma}


Since we have stronger a priori estimates than those in \cite{Alber06} where
 it is assumed that the initial data is in $H^1(\Omega)$, we don't need any more  a compactness lemma of Aubin-Lions type
or its generalized version (see e.g. Simon \cite{Simon} and  Roub\'ic\v ek
\cite{Roubicek}) which   in \cite{Alber06} plays a crucial role in the proof of the existence of weak solutions with $H^{1}(\Omega)$
initial data. A proof of Lemma~\ref{L4.2} can be found e.g.  in the book by Lions
\cite[p.12]{Lions}.

\medskip
\noindent{\it Proof of Lemma~\ref{L4.1}:} Since the estimates in Lemma~\ref{L3.1} and Corollary~\ref{C3.5} imply that the sequence
 $\int_0^{S^\kappa_x}|y|^\frac12_\kappa dy$ is uniformly bounded in $H^1(Q_{T_e})$ for $\kappa\in (0,1]$.
 By the Sobolev imbedding theorem, we assert that $\int_0^{S^\kappa_x}|y|^\frac12_\kappa dy$ is compact
 in $L^2(Q_{T_e}) =  L^2(0,T_e;L^2(\Omega)) $. Thus there is a subsequence, still
denoted by $\int_0^{S^\kappa_x}|y|^\frac12_\kappa dy $, which converges strongly
in $ L^2(Q_{T_e})  $ to a limit function
$G\in L^2(Q_{T_e})  $.  Next we prove that the sequence $\int_0^{S^\kappa_x}|y|^\frac12  dy$
converges to $G$   in $L^2(Q_{T_e}) $. Write
$$
 \int_0^{S_x^\kappa}  |y|^\frac12 dy = \int_0^{S_x^\kappa}  |y|^\frac12_\kappa dy +
 \int_0^{S_x^\kappa}  \Big(|y|^\frac12  - |y|_\kappa^\frac12\Big) dy=I_1+I_2.
$$
It is easy to compute that
\begin{eqnarray}
 0\le |y|^\frac12_\kappa  - |y|^\frac12 &=& \frac{ |y| _\kappa  - |y| }{ |y|^\frac12_\kappa  + |y|^\frac12}\non\\
 &=& \frac{\kappa^2}{\Big(|y|^\frac12_\kappa  + |y|^\frac12\Big)(|y| _\kappa + |y| )}\non\\
 &\le & \frac{\kappa^2}{ \kappa ^{\frac12+1} } = \kappa ^\frac12.
\end{eqnarray}
Thus $I_2$ can be estimated  as
$$
 \|I_2 \|_{L^2(Q_{T_e})} \le  \|\kappa^\frac12 S_x^\kappa\|_{L^2(Q_{T_e})}
 \le C \kappa^\frac12 \|S_{x}^\kappa\|_{L^\infty(0,{T_e};L^2(\Omega))}\le C
 \kappa^\frac12 \to 0 .
$$
Therefore, $\int_0^{S_x^\kappa}  |y|^\frac12  dy\to \lim_{\kappa\to 0} I_1 = G$
strongly in $L^2(Q_{T_e})$.

Consequently, from this  sequence $\int_0^{S_x^\kappa}  |y|^\frac12  dy$ we can select another
subsequence, denoted in the same way, which converges almost
everywhere in $Q_{T_e}$. Using that the mapping
$y\mapsto f(y):=\int_0^{y}|y|^\frac12 dy$
has a continuous inverse $f^{-1}:\cR\to \cR$, we infer that also the
sequence $S^\kappa_{x} = f^{-1}\left(\int_0^{S_x^\kappa}  |y|^\frac12 dy \right)$
converges pointwise almost everywhere to $ f^{-1}(G)$ in
$Q_{T_e}$. From the uniqueness of the weak limit we conclude that $
f^{-1}(G)=S_x$ almost everywhere in $Q_{T_e}$. Thus we prove  (\ref{m4.3a}).

To prove (\ref{m4.5}) we note that the estimate $|S_x^\kappa|_\kappa\le |S_x^\kappa| + \kappa $ and the inequality
(\ref{m4.1}) together imply that the sequence  $|S_x^\kappa|_\kappa
$ is uniformly bounded in $L^2(Q_{T_e})$. Thus, (\ref{m4.5})
is a consequence of (\ref{m4.4}) and Lemma~\ref{L4.2}.

\bigskip
\noindent{\it Proof of Theorem~\ref{T1.1}: } Define the functions $u,T$ by inserting $S$ into
 \eq{m2.8a} and \eq{m2.8}, respectively, where $S$ is the limit function of the sequence $S^\kappa$.
We shall prove that $(u,T, S) $ is a weak solution of problem (\ref{Eeq1a1}) -- (\ref{Eeq4a1}).

Recalling \eq{m2.7} we have $S\in L^\infty(Q_{T_e})$. From this relation, from the
definition of $u,\ T$   we immediately see that $u,\ T$  satisfy (\ref{property01}) and (\ref{property01a}),
 respectively. Observe next that $\|S^\kappa\|_{L^\infty(0,T_e;H^1_0(\Omega))}\le C,$ by
Lemma~\ref{L3.1} and Sobolev's embedding theorem. This implies $S\in
L^\infty(0,T_e;H^1_0(\Omega))$, since we can select a subsequence of
$S^\kappa$ which converges weakly to $S$ in this space. Thus, $S$ satisfies (\ref{property02}).

Noting that from (\ref{convolution}) and (\ref{m4.2})
\begin{eqnarray}
 \|\chi_\kappa*S^\kappa-S \|_{{L^2}(Q_{T_e})}&\le &
 \|\chi_\kappa*(S^\kappa-S)\|_{{L^2}(Q_{T_e})}+
 \|(S-\chi_\kappa*S) \|_{{L^2}(Q_{T_e})}\non\\
 & \le & \|(S-\chi_\kappa*S) \|_{{L^2}(Q_{T_e})}
 + \|S^\kappa-S \|_{{L^2}(Q_{T_e})}\to 0,
 \label{convergConvolution}
\end{eqnarray}
for $\kappa\to 0$, we conclude easily that the function   $(u,T)$  defined in this way satisfy weakly
 equation (\ref{Eeq1a1})  -- \eq{Eeq1aa1}. By definition, if the relation
(\ref{definition}) holds, then the proof of the existence of weak solutions is complete.
 To verify (\ref{definition}) we use that by
construction $ S^\kappa $ solves (\ref{m2.2}). Now we multiply equation (\ref{m2.2}) by a test
function $\varphi \in C_0^\infty((-\infty,T_e)\times\Omega)$ and
integrate the resulting equation over $Q_{T_e}$, then obtain
\begin{eqnarray*}
 0 &=& (S^\kappa_t,\varphi)_{Q_{T_e}} +\left(-c\, \nu |S_x^\kappa |_\kappa  S_{xx}^\kappa
 + {\cal F}^\kappa (|S^\kappa _{x}|_\kappa- \kappa),\varphi \right)_{Q_T}\non\\
 &=& -(S^\kappa_0,\varphi(0))_{\Omega}-(S^\kappa,\varphi_t)_{Q_{T_e}}
 +\left(c\, \nu  \int_0^{S^\kappa_x }|y|_\kappa dy,\varphi_{x}\right)_{Q_{T_e}}\non\\
 & & + \left( {\cal F}^\kappa (|S^\kappa _{x}|_\kappa- \kappa),\varphi \right)_{Q_T},
\end{eqnarray*}
where ${\cal F}^\kappa = -c\,(T^\kappa\cdot\bar\varepsilon - \hat\psi'(S^\kappa))$.
Equation (\ref{definition}) follows from this relation if we show that
\begin{eqnarray}
 (S^\kappa_0,\varphi(0))_{\Omega} &\to&
 (S_0,\varphi(0))_{\Omega},
 \label{converg0}
 \\[1ex]
(S^\kappa,\varphi_t)_{Q_{T_e}} &\to& (S,\varphi_t)_{Q_{T_e}},
 \label{converg1}
 \\[1ex]
 \left(\int_0^{S^\kappa_x }|y|_\kappa dy,\varphi_x\right)_{Q_{T_e}}
 &\to& \left(\frac12|S_x|S_x,\varphi_x\right)_{Q_{T_e}},
 \label{converg2}
 \\[0.2cm]
 \left( {\cal F}^\kappa |S^\kappa _{x}|_\kappa , \varphi\right)_{Q_{T_e}}
 &\to& \left( {\cal F} |S_{x}| , \varphi\right)_{Q_{T_e}},
 \label{converg3}\\[0.2cm]
 \left( \kappa {\cal F}^\kappa  , \varphi\right)_{Q_{T_e}}
 &\to&  0 ,
 \label{converg4}
\end{eqnarray}
for $\kappa\to 0$. Now, the relation (\ref{converg0}) follows from
(\ref{approximate}), and the relation (\ref{converg1}) is a
consequence of (\ref{m4.2}). By \eq{m4.4} and \eq{sxsxx2} from which it is easy to get
$\|\int_0^{S^\kappa_x }|y|_\kappa dy\|_{ L^2(Q_{T_e})}\le C$, using again
 Lemma~\ref{L4.1}, one has \eq{converg2}. Convergence \eq{converg4} follows from \eq{5max} easily.

To verify (\ref{converg3}) we note that \eq{convergConvolution},
\eq{5f},   \eq{5max}, and the definition of ${\cal F}^\kappa$ yield
\begin{eqnarray}
 \|{\cal F}^\kappa |S^\kappa _{x}|_\kappa \|_{ L^2(Q_{T_e})}&\le &C,\\[0.2cm]
 {\cal F}^\kappa |S^\kappa _{x}|_\kappa  &\to& {\cal F} |S  _{x}|  ,\ {\rm almost\ everywhere}.
 \label{m2.8a2}
\end{eqnarray}
Then by Lemma~\ref{L4.1},
$$
 {\cal F}^\kappa |S^\kappa _{x}|_\kappa  \rightharpoonup  {\cal F}  |S_{x}| ,
$$
weakly in $L^2(Q_{T_e})$, which implies (\ref{converg3}). Consequently (\ref{definition}) holds.

\bigskip
\noindent{\bf Regularity.} Since  $S_0\in H^2(\Omega)$, we
can obtain more regular solutions.
By the estimate $\|S^\kappa_t\|_{ L^\infty(0,{T_e};L^2(\Omega))}\le C$, we see that the relation $S_t\in L^\infty(0,{T_e};L^2(\Omega))$ is true.
Then by the  theory of elliptic systems, we obtain (\ref{proper1}).

To prove (\ref{proper2}), we recall the definition of weak solutions. From \eq{definition} it follows that
\begin{eqnarray}
 |( |S_x|S_x,\,\varphi_x )_{Q_{T_e}}   | &\le & C
  \left|\left( \big(T\cdot\overline{\varepsilon} - \hat{\psi}'(S) \big) |S_x|,\,
 \varphi\right)_{Q_{T_e}} \right| + | (S,\varphi_t)_{Q_{T_e}} + (S_0,\varphi(0))_\Omega | \non\\
 &\le & C \|S_x\|_{L^\infty(0,{T_e};L^2(\Omega))} \| \varphi\| _{L^1(0,{T_e}; L^2(\Omega))}  + | (S_t,\varphi )_{Q_{T_e}}   | \non\\
 &\le & C \| \varphi\| _{L^1(0,{T_e}; L^2(\Omega))}  + \| S_t\|_{L^\infty(0,{T_e};L^2(\Omega))} \| \varphi\| _{L^1(0,{T_e}; L^2(\Omega))}\non\\
 &\le & C \| \varphi\| _{L^1(0,{T_e}; L^2(\Omega))},
 \label{definition1}
\end{eqnarray}
here, we used the estimates $\|S_x\|_{L^\infty(0,{T_e};L^2(\Omega))}\le C$ and $\|S_t\|_{L^\infty(0,{T_e};L^2(\Omega))}\le C$.  The right hand side
 of \eq{definition1} holds for all $\varphi\in {L^1(0,{T_e}; L^2(\Omega))}$, whence
\begin{eqnarray}
\sup_{0\le t\le T_e} \|(|S_x|S_x)_x(t)\| &=& \sup_{\|\varphi\|_{L^1(0,{T_e}; L^2(\Omega)) } \le 1} |
 ((|S_x|S_x)_x, \varphi )_{Q_{T_e}}|\non\\
 &=& \sup_{\|\varphi\|_{L^1(0,{T_e}; L^2(\Omega)) } \le 1} |\langle  |S_x|S_x,\,\varphi_x \rangle_{Q_{T_e}} |\non\\
 &\le& C .
 \label{definition2}
\end{eqnarray}
Thus, $(|S_x|S_x)_x\in L^\infty(0,T_e;L^2(\Omega))$.

Furthermore, from the Poincar\'e inequality and the estimate $\|S_x\|_{L^\infty(0,{T_e};L^2(\Omega))}\le C$
 we obtain $|S_x|S_x \in L^\infty(0,T_e; H^{1}(\Omega))$, from which one asserts by the Sobolev imbedding theorem
  that $|S_x|S_x  \in L^\infty(Q_{T_e})$, hence $S_x\in L^\infty(Q_{T_e})$. And the proof of Theorem~\ref{T1.1} is complete.

\bigskip
\bigskip
\noindent{\bf Acknowledgement.} The author of this work has
been partly supported by Grant MTM2008-03541 of
the Ministerio de Educac\'ion y Ciencia of Spain, and by Project
PI2010-04 of the Basque Government.

\end{document}